\documentclass{article} 
\usepackage[utf8]{inputenc}
\usepackage[margin =2.4cm]{geometry}
\usepackage{amssymb}
\usepackage{amsmath}
\usepackage{amsthm}
\usepackage{enumerate} 
\usepackage{graphicx}
\usepackage{csquotes}
\usepackage{subcaption}
\usepackage{authblk}
\usepackage{tikz}
\usepackage{xurl}
\usepackage{hyperref}

\title{\LARGE{When Natural Variables Are Not Enough: Teaching Integer Programming with Sudoku}}

\author{
    {\Large Aled Williams}\\
    Department of Computer Science\\
    Royal Holloway, University of London\\
    London, UK\\
    \texttt{aled.williams@rhul.ac.uk}
}

\begin{document}
	\date{}
	\maketitle

\begin{abstract}
Sudoku is a compact and familiar setting for teaching a surprisingly deep lesson in integer linear programming, namely that the most natural decision variables are not always enough to produce an effective or convenient linear model. This paper compares two formulations of Sudoku. The first uses binary assignment variables indicating whether a particular digit is assigned to a particular cell. This formulation is less natural from the perspective of the puzzle board itself, but it encodes the puzzle rules through simple assignment constraints and extends easily to variants such as Killer Sudoku. The second formulation uses the more natural approach of assigning one integer variable to represent the value in each cell, but then the central requirement is that the values in each row, column, and block must be all different. When this all-different requirement must be expressed using linear constraints while retaining the cell-value variables, the formulation becomes a large collection of pairwise disjunctions linearised by big-$M$ inequalities. The paper includes AMPL code for both the assignment model and the natural all-different model, and derives a second integer program for certifying uniqueness. The paper also uses generalised Sudoku as a careful entry point to computational complexity, while emphasising that the standard nine-by-nine puzzle is not itself an asymptotic problem class.

\vspace{1.0mm}
\end{abstract}


\section{Introduction}
Sudoku has been used as a classroom example in operations research and integer programming because its rules are simple to state, readily illustrated on a familiar puzzle grid, and naturally lead to combinatorial decision models. Chlond~\cite{chlond2005classroom} presents Sudoku as one of two classroom drill exercises in integer programming modelling, using a five-index binary formulation and an Xpress-Mosel implementation. Weiss and Rasmussen~\cite{weiss2007lessons} revisit the puzzle from a spreadsheet-modelling perspective, using a three-index binary assignment formulation and showing how it can be implemented and solved in Microsoft Excel. These contributions show that a familiar recreational puzzle can provide an accessible setting in which students encounter binary variables, assignment constraints, and modelling software.

This paper builds on previous work by making the contrast between natural and effective formulations its main pedagogical focus. In particular, when students first encounter the puzzle, the most literal modelling choice is often to assign one integer variable to each cell, thereby mirroring the puzzle board directly. The difficulty is not in defining these variables, but in expressing the Sudoku rules as linear constraints. The modelling challenge is to ensure that every row, every column, and every block contains distinct values. In other words, the central requirement is an all-different condition, not simply a condition on the sum of the entries. This tension between variables that mirror the puzzle directly and variables that simplify the constraints creates a useful basis for classroom discussion about formulation choice.

This tension highlights an important principle in integer programming, namely that variables that provide the most natural description of a problem are not always those that yield the simplest or most effective linear formulation. The aim of this paper is not merely to show that Sudoku can be solved by integer programming, but to show why variables that describe the solution directly need not lead to the most effective linear formulation.

We first develop the standard binary assignment formulation. We then consider the natural cell-value formulation and examine the linearisation required to model the all-different relation. We also show how the assignment variables support uniqueness certification and Killer Sudoku constraints, and how AMPL can express the natural model through its \texttt{alldiff} operator. Finally, we outline how generalised Sudoku can be used as a careful classroom entry point into $\mathcal{NP}$-completeness.


\section{The Binary Assignment Formulation} \label{sec:binary_assign}
The most widely used integer programming formulation for Sudoku does not model the value written in each cell directly. Instead, it records which digits are assigned to which cells. This formulation is not the most literal representation of the board, but it provides a clean linear representation of the digit-placement decisions underlying the puzzle.

Let $N = \{1,2,\ldots,9\}$, and use this notation throughout the paper. For each cell $(i,j) \in N \times N$ and digit $k \in N$, we define the binary variable
$$
x_{ijk}=
\begin{cases}
1, & \text{if digit } k \text{ is placed in cell } (i,j),\\
0, & \text{otherwise.}
\end{cases}
$$
The model therefore contains $9^3 = 729$ binary variables, each representing a possible assignment of a digit to a cell. 

The Sudoku rules can now be expressed as linear constraints. The requirement that every digit must appear exactly once in each column can be modelled by
$$
\sum_{i\in N} x_{ijk}=1
\qquad \text{for every } j,k\in N.
$$
Similarly, every digit must appear exactly once in each row, which is enforced by
$$
\sum_{j\in N} x_{ijk}=1
\qquad \text{for every } i,k\in N.
$$
Each digit must also appear exactly once in every $3\times 3$ block. Let $a,b\in\{1,4,7\}$ index the upper-left corners of the blocks, then this condition can be modelled via
$$
\sum_{q=0}^{2}\sum_{r=0}^{2} x_{a+q,\,b+r,\,k}=1
\qquad \text{for every } a,b\in\{1,4,7\} \text{ and } k\in N.
$$
Finally, each cell of the Sudoku grid must contain exactly one digit, which yields the constraints
$$
\sum_{k\in N} x_{ijk}=1
\qquad \text{for every } i,j\in N.
$$

Recall that, before incorporating the given clues, the assignment formulation outlined above contains $9^3 = 729$ binary variables and $4\cdot 81=324$ equality constraints, corresponding to the column, row, block, and cell constraints. Although the number of variables is considerably larger than the $81$ cells of the puzzle, the benefit is that every Sudoku rule is represented by a simple linear equality.

Given clues are incorporated naturally. In particular, if the puzzle specifies that cell $(i,j)$ contains the digit $g_{ij} \in N$, then this can be imposed simply by fixing
$x_{ijg_{ij}}=1$.

Observe that Sudoku is fundamentally a feasibility problem, and no objective function is required. In AMPL or other modelling software, one may use a dummy objective such as minimising 0 if needed. 

If desired, the corresponding cell values can be recovered from the binary solution via
$$
z_{ij}=\sum_{k \in N} k \, x_{ijk}, 
$$
where $z_{ij} \in N$ denotes the value in the cell $(i,j)$.
This observation is important pedagogically. The binary variables do not discard the natural cell values, but instead represent each possible assignment of a digit to a cell separately. 
Importantly, although this formulation introduces substantially more variables than there are cells in the puzzle, the resulting model consists entirely of simple linear equalities. This makes the formulation particularly attractive for implementation in standard integer programming software.


\subsection{Implementing the Assignment Model in AMPL} \label{sec:ampl_assign}
The three-index assignment formulation is compact enough to be presented directly to students as an AMPL model. The following implementation uses sets of cell coordinates for each $3\times 3$ block, allowing the block constraints to mirror the row and column constraints.

Before giving the AMPL code, it is useful to note how the mathematical notation given above maps onto the implementation. The set \texttt{VALUES} corresponds to the digit set $N$, while \texttt{ROWS} and \texttt{COLUMNS} index the cell positions. The binary variable \texttt{assign[i,j,v]} is the AMPL counterpart of $x_{ijv}$. The four families of constraints in the algebraic model appear directly as \texttt{row\_sum}, \texttt{col\_sum}, \texttt{block\_sum}, and \texttt{one\_value}. The constraint family \texttt{initial\_value} fixes the variables corresponding to the given clues.

{\small
\begin{verbatim}
set ROWS = 1..9;
set COLUMNS = 1..9;
set VALUES = 1..9;
set BLOCK_ROWS = 1..3;
set BLOCK_COLS = 1..3;
set BLOCK{br in BLOCK_ROWS, bc in BLOCK_COLS} within ROWS cross COLUMNS =
    setof {i in 3*br-2..3*br, j in 3*bc-2..3*bc} (i,j);

param initial{ROWS, COLUMNS} integer >= 0 <= 9 default 0;
var assign{ROWS, COLUMNS, VALUES} binary;

minimize feasibility: 0;

subject to row_sum{i in ROWS, v in VALUES}: 
    sum{j in COLUMNS} assign[i,j,v] = 1;
subject to col_sum{j in COLUMNS, v in VALUES}: 
    sum{i in ROWS} assign[i,j,v] = 1;
subject to block_sum{br in BLOCK_ROWS, bc in BLOCK_COLS, v in VALUES}:
    sum{(i,j) in BLOCK[br,bc]} assign[i,j,v] = 1;
subject to one_value{i in ROWS, j in COLUMNS}: 
    sum{v in VALUES} assign[i,j,v] = 1;
subject to initial_value{i in ROWS, j in COLUMNS, v in VALUES: initial[i,j] = v}:
    assign[i,j,v] = 1;
\end{verbatim}
}

The purpose of giving the AMPL model is not only to provide executable code, but also to show how the formulation can be implemented without losing its mathematical structure. This reflects one of the central motivations for algebraic modelling languages, where they allow a model to be written in a form close to its algebraic statement while remaining precise enough to be solved computationally \cite{fourer1990modeling}. This correspondence is especially clear in our case, as the variable \texttt{assign[i,j,v]} is the AMPL counterpart to $x_{ijv}$, while the row, column, block, and cell rules appear as separate named constraint families whose indexing closely mirrors the mathematical formulation.

This also prepares the comparison with the more natural (or literal) cell-value formulation introduced in the next section. That formulation is often the first one suggested by students, since it assigns a single variable to the value in each cell. The assignment formulation outlined above makes a useful trade-off, replacing each cell by nine binary choices, while turning the Sudoku rules into repeated linear equalities over different index sets. This is pedagogically useful because it makes formulation choice visible, allowing students to compare a variable choice that is natural from the perspective of the puzzle board with one that is less literal but leads to simpler linear constraints. The same assignment formulation is also useful for later extensions. In particular, the uniqueness formulation in Section \ref{sec:uniqueness} is transparent because a completed grid is described by the 81 digit-cell assignments whose variables take value 1, while the Killer Sudoku extension in Section \ref{sec:killer} can express cage sums directly through the same binary variables. 


\section{The Natural Variable Formulation} \label{sec:natural}
We now introduce a more natural (or literal) formulation, often the first suggested by students. Instead of introducing a binary variable for each possible digit-cell assignment, we assign one integer variable to each cell. In particular, let $z_{ij} \in N$ denote the value in cell $(i,j) \in N \times N$. The appeal of this formulation is that it mirrors the Sudoku grid directly, with 81 variables, one for each cell. This makes the formulation useful pedagogically, since students often choose variables before fully considering the implications for the rest of the model. The difficulty is not in defining the variables, but in expressing the Sudoku rules as linear constraints.

A first attempt, often suggested by students, is to require the values in each row, column, and block to sum to 45, since $1+\cdots+9=45$. This gives the linear constraints 
$$
\begin{aligned}
\sum_{j\in N} z_{ij}&=45 \qquad \text{for every } i\in N, \\
\sum_{i\in N} z_{ij}&=45 \qquad \text{for every } j\in N, \\
\sum_{q=0}^{2}\sum_{r=0}^{2} z_{a+q,b+r}&=45 \qquad \text{for every } a,b\in\{1,4,7\}.
\end{aligned}
$$
These constraints are necessary, but they are not sufficient. The grid in which every cell is assigned the value 5 would satisfy all of the sum constraints above, even though it is not a valid Sudoku grid. The actual requirement is stronger, where we need to enforce that the variables in each row, each column, and each block must take all different values. 

Consider the first row, where we require $z_{11},z_{12},\ldots,z_{19}$ to be all different. This all-different condition is equivalent to the 36 pairwise constraints 
$$
|z_{1s}-z_{1t}|\ge 1
\qquad \text{for every } 1\le s<t\le 9.
$$
Observe that each of these constraints is equivalent to the logical disjunction
$$
\text{either } z_{1s}-z_{1t}\ge 1
\quad \text{or}\quad
z_{1t}-z_{1s}\ge 1.
$$
This disjunction can be modelled using big-$M$ constraints. In particular, we introduce an auxiliary binary variable $y^{\mathrm{row}}_{1st}$ for the pair of cells $(1,s)$ and $(1,t)$ and choose a sufficiently large value of $M$ to relax the inactive inequality. Then the disjunction can be expressed as 
$$
\begin{aligned}
z_{1s}-z_{1t} &\ge 1-My^{\mathrm{row}}_{1st}, \\
z_{1t}-z_{1s} &\ge 1-M\left(1-y^{\mathrm{row}}_{1st}\right), \\
y^{\mathrm{row}}_{1st} &\in\{0,1\}, 
\end{aligned}
$$
where $M = 9$ suffices. Observe that if $y^{\mathrm{row}}_{1st} = 0$, then the first inequality enforces $z_{1s}-z_{1t}\ge 1$, while the second inequality is relaxed, since it reduces to $z_{1t}-z_{1s} \ge 1 - 9 = -8$, which necessarily holds for any $z_{1s},z_{1t}\in\{1,\ldots,9\}$. If instead $y^{\mathrm{row}}_{1st} = 1$, then the second inequality enforces $z_{1t}-z_{1s}\ge 1$, while the first inequality is relaxed. Thus, the binary variable $y^{\mathrm{row}}_{1st}$ behaves akin to an on/off switch, selecting which ordering is enforced between the cell values.

This construction must be repeated for every pair of cells in every row, column, and block. A direct (or na\"ive) implementation gives $27\binom{9}{2}=972$ pairwise disjunctions, which follows since there are 9 rows, 9 columns, and 9 blocks, each of which contains $\binom{9}{2}=36$ unordered pairs. 

However, this direct count includes some duplicate pairwise conditions. Note that a duplicate occurs when the same two cells belong to both a row and a block, or to both a column and a block. For example, the cells $(1,1)$ and $(1,2)$ lie in the same row and also in the same $3\times 3$ block, so the condition $z_{11}\neq z_{12}$ is generated twice in the na\"ive count.

We can count these duplicates explicitly. In any fixed row, the nine cells are split across three different blocks, with three cells in each block. Within one such group of three cells, there are $\binom{3}{2}=3$ pairs of cells that share both the row and block. Since there are three such groups in the row, each row contributes $3\binom{3}{2}=9$ row--block duplicates. Thus, across all rows, this gives $9 \cdot 3\binom{3}{2}=81$ row--block duplicates. 

Similarly, each column is also divided among three blocks, with three cells in each block, so there are another 81 column--block duplicates. Thus, after removing these duplicated conditions, the number of distinct pairwise disequalities is $972-81-81=810$. Finally, since each condition is modelled with one auxiliary binary variable and two big-$M$ inequalities, even the formulation after removing duplicated conditions requires 810 auxiliary binary variables and 1620 big-$M$ inequalities.

This is the main lesson of the natural formulation. The variables $z_{ij} \in N$ are natural and few in number, and they match the way the puzzle is usually described. However, the choice of variables must also take into account how the constraints will be expressed. In this case, the central modelling requirement is an all-different relation. If the modelling language or solver interface supports this relation natively, the natural formulation is elegant (as discussed later in Section~\ref{sec:all_diff}). 

If the model must be written as an integer program, then the modeller must either linearise the global constraint, producing many auxiliary binary variables and big-$M$ inequalities, or choose a different variable scheme. This gives students a concrete example of a broader modelling issue, namely that big-$M$ reformulations of logical conditions are often easy to state, but may motivate alternative reformulations when the resulting mixed-integer program is difficult to solve~\cite{codato2006combinatorial}. Thus, the example shows how a formulation based on natural variables may be valid, yet substantially more complex than one based on alternative modelling choices.


\section{Further Extensions}
The advantages of the binary assignment formulation are not limited to the standard Sudoku rules. We now show how the same variables support two natural extensions, namely certifying uniqueness of a completed puzzle and modelling Killer Sudoku cage constraints.

\subsection{Certifying Uniqueness} \label{sec:uniqueness}
Sudoku puzzles are typically published with the promise that the completed grid is unique. The binary assignment formulation allows this property to be verified directly by solving a second integer program.

This can be done more generally for binary integer programs. Suppose we are given a binary integer program of the form
\begin{equation} \label{bip_general}
A \boldsymbol x \le \boldsymbol  b, \qquad \boldsymbol x \in \{0,1\}^n,
\end{equation}
and we are given a feasible solution $\bar{\boldsymbol x}$, i.e. $\bar{\boldsymbol x} \in \{0,1\}^n$ and $A \bar{\boldsymbol x} \le \boldsymbol b$. We wish to determine whether there exists another feasible binary solution $\boldsymbol{x}\in\{0,1\}^n$ with $\boldsymbol{x}\neq\bar{\boldsymbol{x}}$.

A natural way to do this is to maximise the number of entries in which a candidate solution $\boldsymbol x$ differs from the known feasible solution $\bar{\boldsymbol x}$. Since the variables are binary, if $\bar{x}_i = 0$, then disagreement occurs precisely when $x_i = 1$, whereas if $\bar{x}_i = 1$, then disagreement occurs precisely when $x_i = 0$. Thus, the number of differing entries is given by
\begin{equation} \label{disagreement_count}
\sum_{i:\bar{x}_i=0} x_i
+
\sum_{i:\bar{x}_i=1} (1-x_i).
\end{equation}
Note that this expression \eqref{disagreement_count} simply counts the number of positions in which $\boldsymbol x$ differs from $\bar{\boldsymbol x}$. Thus, if the maximum value of this expression over all feasible solutions is zero, then no feasible solution distinct from $\bar{\boldsymbol x}$ exists, and $\bar{\boldsymbol x}$ is the unique solution to the binary integer program \eqref{bip_general}.

The expression \eqref{disagreement_count} can be expanded to yield
$$
\sum_{i:\bar{x}_i=0} x_i
-
\sum_{i:\bar{x}_i=1} x_i
+
\sum_{i:\bar{x}_i=1} 1.
$$
Since $\bar{\boldsymbol x}$ is fixed, the final term is constant and can be omitted without changing the set of optimal solutions. Hence, the objective can equivalently be written as
$$
\max
\left\{
\sum_{i=1}^n d_i x_i
:
Ax \le b,\ x \in \{0,1\}^n
\right\},
$$
where $d_i = 1$ if $\bar{x}_i = 0$, and $d_i = -1$ if $\bar{x}_i = 1$.

We now apply this idea to the Sudoku assignment formulation. We define
$$
d_{ijk}=
\begin{cases}
1, & \text{if } \bar{x}_{ijk}=0,\\
-1, & \text{if } \bar{x}_{ijk}=1.
\end{cases}
$$
We then solve the second integer program 
$$
\max \quad 
\sum_{i \in N}
\sum_{j \in N}
\sum_{k \in N}
d_{ijk} \, x_{ijk}
$$
subject to the Sudoku constraints. Recall that every feasible Sudoku completion contains exactly 81 assignment variables equal to 1, since each of the 81 cells must be assigned exactly one digit. Hence the original solution $\bar{\boldsymbol x}$ attains objective value 
$$
\sum_{i \in N}
\sum_{j \in N}
\sum_{k \in N} 
d_{ijk} \, \bar{x}_{ijk} = -81.
$$
It follows that the Sudoku instance is uniquely solvable if and only if the optimal value of the second problem is exactly -81. In particular, observe that any alternative completion must agree with $\bar{\boldsymbol{x}}$ in fewer than all 81 selected digit-cell assignments, and would therefore obtain an objective value strictly greater than $-81$.


\subsection{Killer Sudoku} \label{sec:killer}
The binary assignment formulation naturally extends to \textit{Killer Sudoku}, a variant in which some groups of cells (often called cages) are labelled with the sum of the digits in those cells. Let $B_1,\ldots,B_m\subseteq N\times N$ denote the cages whose sums are given, and let $s_r$ be the required sum for cage $B_r$. Recall that the value in cell $(i,j)$ is recovered via $\sum_{k\in N} k \, x_{ijk}$, and thus the cage condition is modelled via 
$$
\sum_{(i,j)\in B_r}\sum_{k\in N} k \, x_{ijk}=s_r
\qquad \text{for every } r=1,\ldots,m,
$$
which ensures that the digits assigned to the cells within cage $B_r$ sum to the prescribed total $s_r$. 

We distinguish the cage-sum rule from the common Killer Sudoku convention that no digit may repeat within a cage. If the no-repeat convention is intended and is not already implied by the row, column, and block constraints, it can be modelled via 
$$
\sum_{(i,j)\in B_r} x_{ijk}\le 1
\qquad \text{for every } r=1,\ldots,m \text{ and } k\in N,
$$
which ensures that each digit $k$ appears at most once within cage $B_r$.


\section{AMPL and \texttt{alldiff}} \label{sec:all_diff}
Recall that Section~\ref{sec:natural} showed that the natural formulation, often suggested by students, is easy to describe but difficult to write as an explicit integer program. AMPL provides an alternative way to state this formulation through an all-different operator. The all-different relation is a standard global constraint in constraint programming (see e.g. \cite{rossi2006handbook}), with classic filtering methods developed by R{\'e}gin~\cite{regin1994filtering}. Fourer and Gay~\cite{fourer2002extending} describe extensions of AMPL to support such constraint-programming constructs, including an \texttt{alldiff} operator. The current AMPL MP documentation states that \texttt{alldiff} is a constraint-valued operator that requires all expressions in an indexed collection to take different values \cite{amplmp2026expressions}. Notice that each cell value lies in $\{1,\ldots,9\}$, so requiring the nine values in a row, column, or block to be all different ensures that each digit occurs exactly once.

The integer variable \texttt{value[i,j]} is the AMPL counterpart of $z_{ij}$. The three \texttt{alldiff} constraint families correspond directly to the row, column, and block requirements, while \texttt{natural\_initial\_value} fixes the given clues. The natural formulation can be written in AMPL as follows. 

{\small
\begin{verbatim}
set ROWS = 1..9;
set COLUMNS = 1..9;
set BLOCK_ROWS = 1..3;
set BLOCK_COLS = 1..3;

set BLOCK{br in BLOCK_ROWS, bc in BLOCK_COLS} within ROWS cross COLUMNS =
    setof {i in 3*br-2..3*br, j in 3*bc-2..3*bc} (i,j);

param initial{ROWS, COLUMNS} integer >= 0 <= 9 default 0;
var value{ROWS, COLUMNS} integer >= 1 <= 9;

minimize natural_feasibility: 0;

subject to row_alldiff{i in ROWS}: 
    alldiff {j in COLUMNS} value[i,j];
subject to col_alldiff{j in COLUMNS}: 
    alldiff {i in ROWS} value[i,j];
subject to block_alldiff{br in BLOCK_ROWS, bc in BLOCK_COLS}: 
    alldiff {(i,j) in BLOCK[br,bc]} value[i,j];
subject to natural_initial_value{i in ROWS, j in COLUMNS: initial[i,j] > 0}: 
    value[i,j] = initial[i,j];
\end{verbatim}
}

The code above is a complete AMPL implementation of the natural cell-value model. It has one integer variable for each cell and imposes the row, column, and block rules directly through \texttt{alldiff}. In contrast to the big-$M$ formulation of Section~\ref{sec:natural}, the all-different conditions do not appear as a collection of pairwise disjunctions, and no auxiliary binary variables are required in the AMPL code. This makes the model useful pedagogically, since students can compare the natural formulation as they first describe it with the explicit linearisation required to express that formulation as an integer program.

It is important to note that the \texttt{alldiff} operator is a high-level modelling construct rather than an explicit collection of linear constraints. The current documentation for the AMPL MP solver interface, which translates AMPL models into solver-specific representations, states that no supported solver accepts \texttt{alldiff} natively \cite{amplmp2026expressions}. Instead, the interface automatically reformulates the \texttt{alldiff} constraint into a representation involving simpler constraints before passing the resulting model to the solver. Thus, the concise appearance of the AMPL model should not be interpreted as implying that the resulting solver model is equally concise, since the reformulation may involve additional variables and constraints.

This provides a useful pedagogical comparison. The binary assignment formulation of Section~\ref{sec:binary_assign} is a linear program constructed directly by the modeller. The big-$M$ formulation of Section~\ref{sec:natural} is a manual linearisation of the desired all-different relation. The AMPL \texttt{alldiff} formulation above allows the modeller to state the natural all-different relation directly, while leaving its subsequent reformulation to the modelling system. This comparison helps students distinguish between the mathematical relation being modelled and the formulation ultimately used for computation.


\section{Sudoku as a Gateway to NP-Completeness}
Although the focus of this paper is formulation rather than complexity theory, Sudoku also provides a useful classroom entry point into $\mathcal{NP}$-completeness. Sudoku initially appears to be a recreational puzzle rather than an optimisation problem. However, once phrased as a decision problem, namely whether a given partially filled grid admits a feasible completion, it exhibits several features that are central to integer programming and combinatorial optimisation more broadly.

The first observation is that a proposed completed Sudoku grid is easy to verify. In particular, one simply checks that each row, column, and block contains each symbol exactly once, and that the given clues have been respected. In contrast, determining whether a partially filled grid is completable may require exploring many possible assignments. This gives students a concrete and familiar example of the difference between checking a proposed solution and finding one.

A second observation is that the number of candidate assignments is very large even for the standard puzzle. In particular, if one ignores both the given clues and the Sudoku rules, then a standard $9\times 9$ grid has $9^{81}$ possible assignments of digits to cells. The Sudoku rules eliminate most of these assignments, but the feasibility problem still involves many interacting discrete choices.

There is, however, an important caveat. The standard Sudoku puzzle is a fixed $9\times 9$ problem with symbols $1,\ldots,9$ and nine $3\times 3$ blocks. Since the problem size is fixed, the usual asymptotic language of $\mathcal{NP}$-completeness does not apply directly to this finite class of instances. The correct complexity statement concerns a generalised family of Sudoku problems, namely those with an $n^2\times n^2$ grid, $n\times n$ blocks, and symbols $1,\ldots,n^2$. 
For this generalised family, the Sudoku completion problem is $\mathcal{NP}$-complete \cite{yato2003complexity}. More precisely, Yato and Seta \cite{yato2003complexity} show that the associated search problem is $\mathcal{ASP}$-complete using a solution-preserving polynomial-time reduction from partial Latin square completion, a closely related problem that is itself $\mathcal{NP}$-complete \cite{colbourn1984complexity}.

This distinction between the fixed $9\times 9$ problem and the generalised family is pedagogically useful. The standard $9\times 9$ puzzle is small enough to be modelled and solved as a classroom example, while the generalised version shows how the same modelling structure extends to an $\mathcal{NP}$-complete family of feasibility problems. Sudoku therefore provides a natural connection between integer programming modelling and computational complexity. With this fixed-size caveat made explicit, it also provides an accessible bridge from modelling practice to the intuition behind $\mathcal{NP}$-completeness.


\section{Conclusion}
Sudoku provides an accessible setting in which students can encounter one of the central ideas of integer programming formulation, namely that the variables that most naturally describe a problem are not always those that yield the most effective linear formulation. The cell-value variables mirror the puzzle board directly, but the associated all-different requirements are not straightforward to express with linear constraints without introducing auxiliary variables. The binary assignment formulation, in contrast, introduces more variables and is less literal from the perspective of the puzzle board, but transforms the Sudoku rules into a collection of simple linear equalities.

Comparing these formulations is pedagogically useful because it allows students to see that modelling is not simply the process of translating a verbal description into mathematical notation. Instead, modelling involves making choices about how problem structure should be represented so that constraints can be expressed clearly and handled effectively. The contrast between the assignment model, the explicit big-$M$ linearisation, and the high-level \texttt{alldiff} representation via AMPL allows students to compare the consequences of different modelling choices directly. The accompanying extensions on uniqueness certification and Killer Sudoku show how a well-chosen formulation can be adapted rather than rebuilt.

More generally, the pedagogical value of Sudoku extends beyond the puzzle itself. In particular, students can begin with an intuitive cell-value representation, identify why it does not yield a simple linear formulation, and then compare it with an assignment representation that does. This progression helps students see that modelling is an iterative process in which, before committing to variables, they should consider how those variables will affect the constraints and the computational structure of the problem. This formulation principle applies broadly in operations research, where formulation choices can determine whether a mixed-integer program is compact or effective enough to solve in practice \cite{vielma2015mixed}.


\bibliographystyle{plainurl}
\bibliography{references}

\end{document}